\documentclass[12 pt]{article}

\usepackage{amsmath,amsthm,amssymb,verbatim,graphicx,mathrsfs,fullpage,enumerate,colonequals}

\newtheorem{thm}{Theorem}

\newtheorem{ex}{Example}
\newtheorem{lemma}{Lemma}
\newtheorem{propn}{Proposition}

\newcommand{\R}{\mathbb{R}}
\newcommand{\Z}{\mathbb{Z}}

\newcommand{\N}{\mathbb{N}}
\newcommand{\Q}{\mathbb{Q}}
\newcommand{\dom}{\mathrm{dom}}
\newcommand{\ran}{\mathrm{ran}}

\title{Two Topological Uniqueness Theorems for \\ Spaces of Real Numbers}
\author{
Michael Francis\footnote{Supported by a University of Victoria graduate award} \\
Department of Mathematics and Statistics \\
University of Victoria, P.O. Box 3060 STN CSC \\
Victoria, BC, \textsc{Canada} V8W 3R4
\texttt{mfranc@uvic.ca}
}
\date{November 17, 2011}

\begin{document}

\maketitle

\abstract{\noindent A 1910 theorem of Brouwer characterizes the Cantor set as the unique totally disconnected, compact metric space without isolated points. A 1920 theorem of Sierpi\'{n}ski characterizes the rationals as the unique countable metric space without isolated points. The purpose of this exposition is to give an accessible overview of this celebrated pair of uniqueness results. It is illuminating to treat the problems simultaneously because of commonalities in their proofs. Some of the more counterintuitive implications of these results are explored through examples. Additionally, near-examples are provided which thwart various attempts to relax hypotheses. 
\\

\\

}

\setlength{\parindent}{0pt}
\setlength{\parskip}{1ex plus 0.5ex minus 0.2ex}

\section{Introduction}

The problem of characterizing spaces of real numbers topologically is an old and productive one.  In 1928, Alexandroff \& Urysohn characterized the irrationals as the unique separable, completely metrizable, zero-dimensional space for which every compact subset has empty interior \cite{alexandroff/urysohn}. In 1936, Ward characterized the real line as the unique connected, locally connected, separable, metrizable space for which the removal of any point results in precisely two connected components \cite{ward}. In 1970, Franklin \& Krishnarao improved Ward's result by weakening \emph{metrizable} to \emph{regular}  \cite{franklin/krishnarao}, thereby removing implicit mention of the reals from the characterization. As stated in the abstract, our interest is in two other results from this family. To wit, Brouwer's characterization of the Cantor set \cite{cantor}, and Sierpi\'{n}ski's characterization of the rationals \cite{sierpinski}. Our proof of Sieri\'{n}ski's theorem is modeled on one given in \cite{dasgupta}.

First, we standardize some notation and terminology. Except where otherwise stated, all sets of real numbers are assumed to carry the subspace topology inherited from the usual topology on $\R$. We denote by $C$ the standard ``middle thirds'' set of Cantor, and by $E$ the countable, dense subset of $C$ consisting of endpoints of intervals deleted during the construction of $C$, along with $0$ and $1$. Denote by $\N$ the set of positive integers and by $\mathrm{2}$ the (discrete) 2-point space $\{0,1\}$. It is standard that the product topology on the space $\mathrm{2}^\N$ of infinite binary sequences is homeomorphic to $C$ (the natural homeomorphism is the one sending $(b_i) \in \mathrm{2}^\N$ to  $\sum_{i=1}^\infty \frac{2 b_i}{3^i} \in C$). Since the 2-point discrete space has a natural topological group structure (that of $\Z/2\Z$) the identification $C \cong \mathrm{2}^\N$ shows that $C$ has the structure of an uncountable (abelian) topological  group. In particular, this shows $C$ has many self homeomorphisms - a point which will later be of use. The separation axioms \emph{regular} and \emph{normal} are taken to include the Hausdorff condition, by definition. The word \emph{countable} is taken to mean countably infinite.

The overarching goal is to prove the following theorems which characterize $C$ and $\Q$ uniquely up to homeomorphism.

\begin{thm}[Brouwer] Every nonempty, totally disconnected, compact, metrizable space without isolated points is homeomorphic to the Cantor set $C$.
\end{thm}

\begin{thm}[Sierpi\'{n}ski] Every countable, metrizable space without isolated points is homeomorphic to the rational numbers $\Q$.
\end{thm}

Some applications of these results are given in Examples~\ref{brouw ex} and~\ref{sier ex} below.

\begin{ex}\label{brouw ex} \text{}
\begin{enumerate}[(a)] 
\item Any nonempty perfect subset of $C$ is homeomorphic to $C$.
\item Extending (a), a nonempty subset of $\R$ is homeomorphic to $C$ if and only if it is perfect, nowhere dense and compact.
\item The products spaces $C^2,C^3, \ldots$ and $C^\N$ are homeomorphic to $C$.
\item If $X_1,X_2,\ldots$ are totally disconnected, compact, metrizable spaces (e.g. finite discrete spaces) with at least 2 points, then $\prod_{i =1}^\infty X_i$ is homeomorphic to $C$.
\end{enumerate}
\end{ex}

Part (d) holds because total disconnectedness and compactness are preserved under arbitrary products, metrizability is preserved under countable products, and the product of an infinite family of spaces with at least 2 points has no isolated point.

\begin{ex}\label{sier ex} \text{}
\begin{enumerate}[(a)]
\item Any countable, dense subset of $\R$ is homeomorphic to $\Q$ (already for spaces as simple as $\Q \setminus \{0\}$ or $\Q \cup \{\sqrt 2\}$ this is not completely obvious). 
\item The endpoints $E$ of the Cantor set are homeomorphic to $\Q$.
\item Extending (a), a countable dense subset of any Euclidean $n$-space is homeomorphic to $\Q$. In particular, the product spaces $\Q^2,\Q^3,\ldots$ are homeomorphic to $\Q$.
\item The ``Sorgenfrey'' topology on $\Q$ generated by half-open intervals $[a,b) \cap  \Q$ is homeomorphic to the standard one (metrizability is most easily seen through Urysohn's metrization theorem).
\end{enumerate}
\end{ex}

The spaces in Part (a) above are more basic in the sense that it will in fact be possible to find an \emph{order preserving} homeomorphism with $\Q$ for these spaces. In Part (b), this cannot be so because the order type of $E$ is very different from that of $\Q$ (every left endpoint is the immediate predecessor of a right endpoint). Thus, an order preserving homeomorphism is too much to hope for in general. This concession paves the way for the intuition-testing examples in Part (c). Here, it is far from apparent that there should be any order structure whatsoever compatible with the topology. Part (d) may be the most perverse example of the lot; the Sorgenfrey topology on $\Q$ is strictly finer topology than the standard topology while, at the same time, homeomorphic to the standard topology.

An order preserving homeomorphism between two subsets of $\R$ will be called an \emph{order homeomorphism}. Our general outline for proving Brouwer's theorem and Sierpi\'{n}ski's theorem is the same:

\begin{enumerate}
\item Identify the subsets of $\R$ that are order-homeomorphic to $C$ and $\Q$ respectively.
\item Given $X$ as in Brouwer's theorem or Sierpi\'{n}ski's theorem, construct an embedding of $X$ into $\R$ whose range is one of the above sets.
\end{enumerate}

\section{The Ordered Perspective}

In this paper, we commit a mild travesty by only considering suborders of $\R$. When $X \subset \R$, there is a second natural topology on $X$, besides the subspace topology. This is the \emph{order topology} on $X$. The sets
\[ (a,b) \cap X \ \ \  [\ell,b) \cap X \ \ \ (a,r] \cap X \]
where $\ell$ and $r$ are the largest and smallest elements of $X$ (if indeed such exist) and $a,b \in X$ satisfy $\ell \leq a < b \leq r$ are a basis for the order topology.  Note that an order isomorphism between $X,Y \subset \R$ is also a homeomorphism of their order topologies since basic sets will be sent to basic sets. The subspace topology on $X \subset \R$ is always finer than the order topology and, in many familiar instances, the distinction is irrelevant. This is the case for $C$ and $\Q$ (if the former is not clear now, it will be made clear soon). It is an inconvenient fact of life that the subspace topology can be strictly finer. For example, there is a clear order isomorphism between the spaces

\[ X = [0,1) \cup \{2\} \cup (3,4] \text{ and } Y = [0,2] \]

so their order topologies are homeomorphic. However, their subspace topologies are not homeomorphic since $X$ has an isolated point and $Y$ does not. We will record an easy, but highly checkable, necessary and sufficient condition for these topologies to agree. It seems useful and appropriate to phrase this condition, and two others that will become relevant, in a common language. 

By a \emph{gap} in $X$ we mean a bounded connected component of $\R \setminus X$. We refine this notion. Every bounded, connected subset of $\R$ has precisely one of the following forms
\[ \underbrace{(a,b)}_\text{open} \ \ \ 
\underbrace{[a,b] \ \ \{a\}}_\text{closed} \ \ \ 
\underbrace{[a,b) \ \ (a,b]}_\text{half-open}, \]
so we may classify gaps disjointly as follows: an open gap is an \emph{essential gap} ($X$ contains both endpoints), a closed gap is a \emph{Dedekind gap} ($X$ contains neither endpoint), a half-open gap is a \emph{pseudo-gap} ($X$ contains one endpoint and not the other). For example 
\[ X = ([0,1] \cap \Q) \cup [2,3] \cup (4,5) \cup(5,6] \cup [7,\infty) \]
has 2 essential gaps, many +1 Dedekind gaps and 1 pseudo-gap. The first two flavors of gaps are intrinsic to the order on $X$. The essential gaps in $X$ are in one to one correspondence with pairs $(x,y)$ of points in $X$ where $x$ is an immediate predecessor of $y$ (so $y$ is the immediate successor of $x$).  The Dedekind gaps in $X$ are in one to one correspondence with the \emph{Dedekind cuts} in $X$, where a Dedekind cut in $X$ is defined as a pair $(L,U)$ of nonempty subsets of $X$ such that $L \cup U =X$, $L<U$, $L$ has no largest element, and $U$ has no smallest element. In contrast, however, pseudo-gaps, are \emph{not} detected by the order on $X$ and, when present, indicate that the way $X$ sits in $\R$ is somehow defective. Their importance stems from the following easy result.

\begin{propn}\label{top agree}
The order topology on $X \subset \R$ agrees with the subspace topology if and only if $X$ has no pseudo-gaps.
\end{propn}

Note that a closed subset of $\R$ can only have essential gaps (it must contain both endpoints of any gap). In particular, this shows the order topology and subspace topologies on $C$ are the same. 

Essential gaps, on the other hand, are important to us because of their role they play in the following classical, order theoretic analog of Sierpi\'{n}ski's theorem - originally due to Cantor \cite{cantor}.

\begin{thm}[Cantor]\label{cantor}
If $A,B \subset \R$ are countable, have no essential gaps\footnote{It is standard to call a linear order with no essential gaps \emph{dense}, but this conflicts with the meaning of dense in topology so we avoid this usage here.}, and have neither largest nor smallest elements (these conditions mean, for instance that, when $a,a' \in A$, $a<a'$, we can find $x,y,z \in A$ such that $x<a<y<a'<z$), then $A,B$ are order isomorphic. In particular, both are order isomorphic to $\Q$.
\end{thm}
\begin{proof}
First we fix enumerations of $A$ and $B$. It simplifies matters to index the elements of $A$ with odd numbers and the elements of $B$ with even numbers. So, $A = \{a_1,a_3,a_5,\ldots\}$ while $B = \{b_2,b_4,b_6,\ldots\}$. Let $\mathscr{F}$ be the collection of partially defined, order preserving maps $A \to B$ with finite domain. The plan is to build up an order isomorphism through a countable sequence of extensions $\varphi_1,\varphi_2,\varphi_3,\ldots \in \mathscr{F}$. Moreover, we choose our extensions so that 
\[ a_1 \in \dom \varphi_1 \ \ b_2 \in \ran \varphi_2 \ \ a_3 \in \dom \varphi_3 \ \ b_4 \in \ran \varphi_4 \ \ \ldots \]
which ensures the map defined ``in the limit'' is defined on all of $A$ and hits all of $B$. This is sometimes called a ``back-and-forth'' construction. This process can always continue because of the following.

\emph{Claim.} If $\varphi \in \mathscr{F}$ and $a \in A$, then there is an order preserving extension $\varphi^*$ of $\varphi$ with $\dom \varphi^* = \dom \varphi \cup \{a\}$.  If $\varphi \in \mathscr{F}$ and $b \in B$, then there is an order preserving extension $\varphi^*$ of $\varphi$ with $\ran \varphi^* = \ran \varphi \cup \{b\}$.

We prove only the first statement. If $a \in \dom \varphi$ already, then do nothing. If $a < \dom \varphi$, then we use the assumption that $B$ has no smallest element to produce $b \in B$ with $b <\ran \varphi$ (note $\ran \varphi$ is finite) and define $\varphi^*(a) = b$. We proceed similarly when $a > \dom \varphi$. Otherwise, $a$ has an immediate predecessor $a^-$ and an immediate successor $a^+$ in $\dom \varphi$ (since $\dom \varphi$ is finite). Since $\varphi(a^-) < \varphi(a^+)$ and $B$ has no essential gaps, there is a $b \in B$ with $\varphi(a^-) < b < \varphi(a^+)$ and we set $\varphi^*(x) = b$. In all of the above cases, the map $\varphi^*$ so obtained is an extension of the desired type.
\end{proof}

We now identify the subsets of $\R$ which are order-homeomorphic to $\Q$ (i.e. homeomorphic to $\Q$ via. an order preserving homeomorphism). Let $X \subset \R$. We say $p \in \R$ is a \emph{left limit point} of $X$, if it is a limit point of $(-\infty,p) \cap X$. A \emph{right limit point} of $X$ is defined similarly. If $p$ is both a left limit point and a right limit point of $X$, then we say that $p$ is a \emph{2-sided limit point} of $X$.

\begin{thm}[Ordered version of Sierpi\'{n}ski's theorem]\label{ord sier}
A set $X \subset \R$ is order-homeomorphic to $\Q$ if and only if $X$ is countable and every $x \in X$ is a 2-sided limit point of $X$. 
\end{thm}
\begin{proof}
Only the ``if'' part of the statement is nontrivial. Suppose $X$ is countable and each point of $X$ is a 2-sided limit point of $X$. No $x \in X$ can be the left endpoint of a gap in $X$ or a largest element of $X$, or else $x$ is not a right limit point of $X$. Similarly, no $x \in X$ can be a the right endpoint of a gap in $X$ or a smallest element of $X$, or else $x$ is not a left limit point if $X$. It follows that $X$ has no largest or smallest element, no essential gaps, and no pseudo-gaps. By Cantor's theorem, there is an order isomorphism $\varphi : X \to \Q$ which is homeomorphism from the the order topology on $X$ to $\Q$. However, since $X$ has no pseudo-gaps, Theorem~\ref{top agree} shows $\varphi$ is also a homeomorphism from the subspace topology on $X$ to $\Q$.
\end{proof}

Before we can prove an analogous result for $C$, we will need to our attention to Dedekind gaps. We say that $X \subset \R$ is \emph{Dedekind complete} if it has no Dedekind gaps. As was previously observed, a closed subset of $\R$ can have only essential gaps. Thus, closed subsets of $\R$ are Dedekind complete. It turns out that any linear order embeds into a certain essentially unique Dedekind complete linear order which could be called its \emph{Dedekind completion} (the construction just mimics the construction of the reals using Dedekind cuts in $\Q$). We stop short of defining a Dedekind completion precisely, or proving it always exists. However, we will use the following result - a special case of uniqueness for Dedekind completions.

\begin{propn}\label{ded compl}
Suppose $A \subset X \subset \R$, $X$ is Dedekind complete, and each point in $X \setminus A$ is a 2-sided limit point of $A$. Suppose $B \subset Y \subset \R$, $Y$ is Dedekind complete, and each point in $Y \setminus B$ is a 2-sided limit point of $B$. Then any order isomorphism $\varphi : A \to B$ has a unique extension to an order isomorphism $\varphi^* : X \to Y$. 
\end{propn}
\begin{proof}
Let $x \in X \setminus A$. Let $L= \{y \in Y : y \leq \varphi(a) \text{ for some } a \in A \text{ with } a < x\}$ and $U = \{y \in Y : y \geq \varphi(a) \text{ for some } a \in A \text{ with } a > x\}$. Since $x$ is a 2-sided limit point of $A$ and since $\varphi$ is order preserving,  it follows that $L$ and $U$ are nonempty, $L < U$, $L$ has no largest element, and $U$ has no smallest element. Since $Y$ is Dedekind complete, $(L,U)$ cannot be a Dedekind cut in $Y$ and there exists $y \in Y$ satisfying $L < y< U$. In fact, the $y$ must be unique since if $y_1,y_2$ have $L <  y_1  <  y_2  < U$, then $y_2$ is not a left limit point of $B$. Clearly we are forced to define $\varphi^*(x) = y$ if the extension is to be order preserving. The map $\varphi^*$ obtained by carrying out this argument at each $x \in X \setminus A$ is an order preserving injection $\varphi^* : X \to Y$ extending $\varphi$. The proof that $\varphi^*$ is onto proceeds similarly by fixing a $y \in Y \setminus B$. 
\end{proof}

Now we identify the subsets of $\R$ which are order-homeomorphic to $C$. 

\begin{thm}[Ordered version of Brouwer's theorem]\label{ord brouw}
A nonempty set $X \subset \R$ is order-homeomorphic to $C$ if and only if $X$ is perfect, nowhere dense and compact.
\end{thm}
\begin{proof}
Only the ``if'' part of the statement is nontrivial. Let $X$ be as above. Since $X$ is closed, it has only essential gaps. Let $L, R \subset X$ be, respectively, the set of left endpoints and the set of right endpoints of essential gaps in $X$. Note that $L \cap R = \varnothing$ since $X$ has no isolated points.  

First we show  $L$ is order-homeomorphic to $\Q$ by applying Cantor's theorem. Consider the left endpoint $\ell$ of some essential gap $(\ell, r)$ in $X$. It must be that $\ell$ is a left limit of $X$ (or else it is isolated, but $X$ is perfect). Whenever $x \in X$ and $x  < \ell$, there must be an essential gap in $X$ between $x$ and $\ell$ (or else $X$ has nonempty interior, contradicting nowhere denseness) so there is a $\ell^- \in L$ with $x < \ell^- < \ell$. From these observations, it follows that $\ell$ is not a smallest element of $L$ and has no immediate successor. Similarly, $r$ is a right limit point of $X$, and, if $x \in X$ and $\ell<r<x$, there is an $\ell^+ \in L$ with $x < r < \ell^+ < x$. From this it follows that $\ell$ is not a largest element of of $L$ and that $\ell$ has no immediate predecessor. From Cantor's theorem, it follows it now follows that $L$ is order isomorphic to $\Q$ as claimed.

There is a natural order isomorphism between $L$ and $R$ by pairing the left and right endpoint of each essential gap. In fact, $L \cup R$ is order isomorphic to $L \times \{0,1\} \cong \Q \times \{0,1\}$ in the dictionary order (adding the elements of $R$ amounts to adjoining an immediate successor to each element of $L$). Let $a = \inf X$ and $b = \sup X$. Since $X$ is compact, $a,b \in X$. Let $E_X = L \cup R \cup \{a,b\} \subset X$. Clearly the order type of $E_X$ is just $\Q \times \{0,1\}$ with a largest and smallest element adjoined.  In particular,  there is an order isomorphism $\varphi : E_X \to E$.

We claim any point $x \in X \setminus E_X$ is a 2-sided limit point of $E_X$. Such an $x$ is at least a 2-sided limit point of $X$ (since it is not the largest element, not the smallest element, and not an endpoint of a gap). Moreover, if $y \in X$ and $y \neq x$, then there is a gap in $X$ between $x$ and $y$ (since $X$ has empty interior). Thus, there are points of $E_X$ between $x$ and $y$, and it follows that $x$ is a 2-sided limit point of $E_X$ too. But now, Proposition~\ref{ded compl} applies, and $\varphi$ extends uniquely to an order isomorphism $\varphi^* : X \to E$. Na\"{i}vely, need just be a homeomorphism of the order topologies on $X$ and $C$ but, since $X,C$ are closed, these are the same as the subspace topologies and we have our desired order-homeomorphism.
\end{proof}

We have now completed step one from our outline in the introduction. As was observed in Example~\ref{brouw ex}, the subsets of $\R$ which satisfy the hypotheses of Brouwer's theorem are precisely the sets appearing in Theorem~\ref{ord brouw} above\footnote{In contrast, there exist subsets of $\R$ (for example, $E$ or $\Q \cap [0,1]$) which satisfy the hypotheses of Sierpi\'{n}ski's theorem, but are not order-homeomorphic to $\Q$.}.  Revisiting the aforementioned outline, it is now apparent that proving Brouwer's theorem is a slightly less delicate business than we originally supposed. \emph{Any} embedding at all that we construct for Step 2 will hit a set of the desired type.

\section{Embeddings of Zero-Dimensional Spaces}

If we are to have any hope of completing Step 2 from our outline, we will need techniques for embedding spaces into $\R$. Typical embedding theorems in general topology take a family $f_i:X \to Y_i$ of continuous functions on a space $X$ and use them as the coordinate functions of a map from $X$ into the product space $\prod_i Y_i$. For this reason, it tends to be easier to construct embeddings when the target space is a large topological product.  Although $\R$ is not itself of this form, it has a subspace of this form\footnote{Actually, there is another famous subspace of $\R$ homeomorphic to a large product available. The set of irrational numbers in $\R$ is homeomorphic to $\N^\N$ (roughly, via continued fraction representations).} - namely $C \cong \mathrm{2}^\N$. From this point of view, looking specifically at embeddings into $C \subset \R$ is quite a natural thing to try. To clarify what follows, we state an amusing and somewhat neglected criterion for map to be an embedding. The trivial proof is left to the reader.

\begin{lemma}
A continuous function $f:X \to Y$ with $X$ a  $T_0$ space is an embedding if and only if $p \notin  \overline S$ implies $f(p) \notin \overline{f(S)}$ for all $p \in X$, $S \subset X$.
\end{lemma}

We say an indexed family $f_i:X \to Y_i$, $i \in I$ of continuous functions \emph{separates points from closed sets} if for every $p \in X$ and every neighbourhood $U$ of $p$, there exists an index $i \in I$ and a closed set $0 \subset Y_i$ such that $f_i(p) \notin 0$ and $f_i(x) \in 0$ for all $x \in X \setminus U$. Intuitively, one interprets the latter condition as saying that $f_i$ is nonzero at $x$ and vanishes outside of $U$. We now state a quite general embedding theorem, whose proof is fairly transparent view of the preceding lemma.

\begin{thm}\label{embed}
If $X$ is a $T_0$ space and $f_i:X \to Y_i$, $i \in I$ is a family of continuous functions which separates points from closed sets, then the function $f : X \to \prod_{i \in I} X_i$ sending $x$ to $(f_i(x))_{i \in I}$ is an embedding.
\end{thm}

If we want use the preceding to embed a space $X$ into $\mathrm{2}^\N$, we will need lots of continuous functions $X \to \mathrm{2}$. Since $\mathrm{2}$ is discrete, a clopen subset $Q$ of a topological space $X$ is essentially the same as a continuous function $X \to \mathrm{2}$. One simply considers the function which is $1$ on $Q$ and $0$ on $X \setminus Q$. In fact, as is easily verified, this correspondence is such that a collection $Q_i$ of clopen sets is a basis if and only if the corresponding family of functions $f_i : X \to \mathrm{2}$ separates points from closed sets. We say that $X$ is \emph{zero-dimensional} when there is exists a basis of clopen sets. We now characterize the subspaces of $C$.

\begin{propn}\label{cantor embed}
A topological space $X$ can be embedded into $C$ if and only it $X$ is $T_0$, 2nd countable, and zero-dimensional.
\end{propn}
\begin{proof}
Only the ``if'' part of the proposition is nontrivial, so suppose $X$ is $T_0$, 2nd countable and zero-dimensional. Let $\mathscr{B}$ be a countable basis for $X$. Construct a countable collection of clopen sets $\mathscr{Q}$ by choosing, for each pair $U,V \in \mathscr{B}$ where this is possible, a clopen set $Q_{U,V}$ satisfying $U \subset Q_{U,V} \subset V$. It is easy to check $\mathscr{Q}$ is a countable basis of clopen sets. So, there is a countable family $f_1,f_2,\ldots$ of continuous functions  $X \to \mathrm{2}$ which separates points from closed sets.  By Theorem~\ref{embed}, we obtain an embedding of $X$ into $\mathrm{2}^\N \cong C$ so we are done. 
\end{proof}

Note that dropping 2nd countability in the above gives us the spaces embeddable into $\mathrm{2}^I$  for some, possibly uncountable, index set $I$. 

\section{The theorems of Brouwer and Sierpi\'{n}ski}

In this section we complete the proofs of the two main results. In both cases, we prove the space under consideration is zero-dimensional and apply Proposition~\ref{cantor embed}. For Brouwer's theorem, this is slightly delicate. Certainly it is worth pointing out at this point that the following statement is generally false: if $x$ and $y$ are in different connected components of a space $X$ then there exists a separation $U,V$ of $X$ with $x \in U$ and $y \in V$.

\begin{ex}
Let $X$ be the subspace of the Euclidean plane consisting of the points $p=\langle 0,0 \rangle$ and  $q = \langle 0,1  \rangle$ together with the vertical lines $\{1/n\} \times [0,1]$, $n \in \N$. The connected components of $X$ are the individual lines and the singletons $\{p\}$ and $\{q\}$. However, there is no separation $U,V$ of $X$ with $p \in U$ and $q \in V$.
\end{ex}

\begin{proof}[Proof of Brouwer's theorem]
Let $X$ be a nonempty, totally disconnected, compact metric space with no isolated points. We need  only prove that $X$  is zero dimensional since then Proposition~\ref{cantor embed} and Theorem~\ref{ord brouw} prove that $X$ is homeomorphic to $C$. In fact, this all comes down to the following claim, and the rest is a routine compactness argument.

\emph{Claim. }If $x,y \in X$ are distinct, there is a separation $U,V$ of $X$ with $x \in U$, $y \in V$.

Fix $x \in X$ and let $Y$ be the (closed) intersection of all clopen neighbourhoods of $x$. We need to prove that $Y = \{x\}$ and, since $X$ is totally disconnected, it suffices to see $Y$ is connected. To this end, suppose that $A,B$ are disjoint closed sets whose union is $Y$. With no harm done, suppose $x \in A$.  We will show that $B = \varnothing$ so that $Y$ is connected. Since $X$ is a normal space, there exist disjoint open sets $U,V$ such that $A \subset U, B \subset V$. It is clear that the clopen sets which \emph{exclude} $x$ are an open cover of the (compact) set $X \setminus (U \cup V)$. It follows that are finitely many clopen neighbourhoods of $x$ whose intersection $Q$, another clopen neighbourhood of $x$, has $A \cup B \subset Q \subset U \cup V$. But now notice that $Q \cap U = Q \setminus V$ is also a clopen neighbourhood of $x$ so $Y = A \cup B \subset Q \subset U$ requiring $B = \varnothing$. 
\end{proof}

For Sierpi\'{n}ski's theorem, it is more straightforward to check zero-dimensionality, but even after applying Proposition~\ref{cantor embed}, more work needs to be done to get an embedding whose range satisfies the hypothesis of Theorem~\ref{ord sier}.

\begin{proof}[Proof of Sierpi\'{n}ski's theorem]
Let $X = (X,d)$ be a countable metric space without isolated points. Let $D \subset [0, \infty)$ equal the set of distances $d(x,y)$ as $x,y$ range over $X$. Since $X$ is countable, $D$ is countable. Therefore, there exist positive real numbers $\epsilon_1,\epsilon_2,\ldots$ in $(0,\infty) \setminus D$ converging to zero. For $x \in X$ and $n \in \N$, let $U(x,n) = \{ y \in X : d(x,y) < \epsilon_n\}$. By design, the complement of $U(x,n)$ is $\{y \in X : d(x,y) > \epsilon_n\}$ so each $U(x,n)$ is clopen. Moreover, the $U(x,n)$ are a basis for $X$ (even a countable one). So, by Proposition~\ref{cantor embed} there exists an embedding $f:X \to C$. We would like to apply Theorem~\ref{ord sier} to $f(X)$, but this is not yet justified since $f(X)$ could equal, say, $E$ and fail to be order isomorphic to $\Q$. 

Note, however, that $\overline{ f(X) }\subset C \subset \R$ satisfies the hypotheses of Brouwer's theorem, or even Theorem~\ref{ord brouw}. It follows that $f(X)$ is homeomorphic to $C$. Therefore, there is in fact an embedding $g: X \to C$ whose image $g(X)$ is \emph{dense} in $C$. Now recall that $C$ is homeomorphic to $\mathrm{2}^\N$ so that $C$ has a natural topological group structure. We claim that there exists an $x \in C$ which ``translates'' $f(X)$ away from the problematic points of $E$. That is, there exists $x \in C$ such that $(f(X) + x) \cap E = \varnothing$. In fact, whenever $A,B \subset C$ are countable, there must be an $x \in C$ with $(A + x) \cap B = \varnothing$. This is because $(A+x) \cap B  \neq \varnothing$ if and only if $x$ is in the countable set $\{b - a : a, \in A, b \in B\}$ so the set of $x \in  C$ that \emph{don't} work is countable. So, by composing $g$ with an appropriate self homeomorphism of $C$, we obtain an embedding $h:X \to C$ such that $g(X)$ is dense in $C$ and $g(X) \cap E = \varnothing$. Now, we claim that every point in $h(X)$ is a 2-sided limit point of $h(X)$ so that Theorem~\ref{ord sier} can be applied. Indeed, any point $c \in C \setminus E$ is a 2-sided limit of $C$, hence a 2-sided limit point of $h(X)$ (since $h(X)$ is dense in $C$). Therefore, $h(X)$ is (order) homeomorphic to $\Q$. Since $X$ is homeomorphic to $h(X)$, this proves Sierpi\'{n}ski's result.
\end{proof}

\section{Modifying hypotheses}

Sierpi\'{n}ski's theorem characterizes $\Q$ by the properties \emph{countable}, \emph{metrizable} and \emph{no isolated points}. If one wishes to expunge any reference to a metric from the theorem, one can replace \emph{metrizable} with \emph{1st countable and regular}. This is possible because, for countable spaces, 1st countable implies 2nd countable so Urysohn's metrization theorem gives the nontrivial half of the equivalence. The following example shows that \emph{countable}, \emph{regular} and \emph{no isolated points} do not suffice. Not even if \emph{regular} is improved to \emph{zero-dimensional}.

\begin{ex}
Let $X = \Q[x]$ be the space of polynomial functions from $\R \to \R$ with rational coefficients. Give $X$ the topology of pointwise convergence. Clearly $X$ is countable and has no isolated points. Also, $X$ is zero dimensional. To see this, suppose that $U$ is a neighbourhood of $p \in X$. Without loss of generality, there is an $\epsilon >0$ and $a_1,\ldots,a_n \in \R$ such that $U = \{ q \in X : |q(a_i) - p(a_i)| < \epsilon  \text{ for } i=1,\ldots,n\}$. For each $i$, $S_i \colonequals \{q(a_i) : q \in X\}$ is countable, so there exist $r_i,s_i \in \R \setminus S_i$ with $p(a_i) - \epsilon < r_i < p(a_i) < s_i < p(a_i) + \epsilon$. Then, the set $Q \colonequals \{q \in X : r_i < q(a_i) < s_i \text{ for } i=1,\ldots,n\}$ is clopen and satisfies $p \in Q \subset U$, so $X$ is zero dimensional. However, $X$ does not have strong enough countability properties to be metrizeable so is not homeomorphic to $\Q$. It is not difficult to see $X$ is not 1st countable. In fact more is true. Sequences do not suffice to detect limit points in $X$. Let $Y$ be the set of $q \in X$ with $|q(x)| \leq 1$ on $[0,1]$ and with $\int_0^1 q(x) \ dx \geq 1/2$. It is clear that the zero polynomial is in the closure of $Y$. But, the zero cannot be a the limit of a sequence $(q_i)$ in $Y$. If such $q_i$ converge pointwise to zero, then  Lebesgue's dominated convergence theorem implies that the sequence of integrals $\int_0^1 q_i(x) \ dx$ converges to zero as well, which is impossible by design.
\end{ex}

The following example shows that \emph{countable}, \emph{2nd countable} and \emph{no isolated points} do not suffice. Not even if we assume \emph{totally disconnected} and \emph{Hausdorff} in addition.

\begin{ex}
As a set, take $X$ to be the union of $\Q$ with two idealized points $p_0$ and $p_1$. For $n=0,1,2,\ldots$, let $I_n = (n,n+1)$. We topologize $X$ by taking $U \subset X$ open if and only if all of the following hold.
\begin{itemize}
\item $U \cap \Q$ is open in the standard topology on $\Q$.
\item If $p_0 \in U$, then $U$ contains all but finitely many of $I_0,I_2,I_4,\ldots$.
\item If $p_1 \in U$, then $U$ contains all but finitely many of $I_1,I_3,I_5,\ldots$. 
\end{itemize}
It is easy to see that $X$ is 2nd countable, Hausdorff, totally disconnected and has  no isolated points. However, $X$ is not regular. For example, $\Z \subset X$ is closed and $p_0 \notin \Z$, but $\Z$ and $p$ cannot have disjoint neighbourhoods.
\end{ex}

Brouwer's theorem characterizes $C$ by the properties \emph{totally disconnected}, \emph{compact}, \emph{metrizable} and \emph{no isolated points}. Here, \emph{metrizable} can be replaced with \emph{Hausdorff and 2nd countable} (since these conditions are equivalent in the presence of compactness) if a ``metric free'' characterization is desired. Obvious examples show no one of these conditions can be dropped. Finally, in the proof of Brouwer's theorem, compactness was used in a key way to deduce zero-dimensionality from total disconnectedness. The compactness assumption is crucial. For example, Cantor's leaky tent is a (noncompact) subspace of the Euclidean plane which is totally disconnected with no isolated points, but not zero-dimensional - nor even totally separated \cite{ctrexamples}.


\begin{thebibliography}{9}

\bibitem{alexandroff/urysohn}
P. Alexandroff \& P. Urysohn, 
\emph{\"{U}ber nulldimensionale Punktmengen}, 
Math. Ann. \textbf{ 98} (1928), 89-106.

\bibitem{brouwer}
 L.E.J. Brouwer, 
\emph{On the structure of perfect sets of points, }
Proc. Akad. Amstersdam \textbf{12}
(1910),
785-794.

\bibitem{cantor}
G. Cantor, 
\emph{Beitr\"{a}ge zur Begr\"{u}ndung der transfiniten Mengenlehre I}, 
Math. Ann. \textbf{46} (1895), 481-512.

\bibitem{dasgupta}
A. Dasgupta, 
\emph{Countable metric spaces without isolated points},
 Topology Atlas
 (2005).

\bibitem{franklin/krishnarao}
S. P. Franklin \& G. V. Krishnarao,
\emph{On the topological characterization of the real line},
J. London Math. Soc(2), \textbf{2} (1970), 589-591.

\bibitem{sierpinski}
W. Sierpi\'{n}ski, 
\emph{Sur une propri\'{e}t\'{e} topologique des ensembles d\'{e}nombrables denses en soi}, 
Fund. Math. \textbf{1}
(1920), 11-16.

\bibitem{ctrexamples}
L. A. Steen \& J. A. Seebach,
\emph{Counterexamples in Topology},
2nd ed., Springer-Verlag, New York (1978) p. 145.

\bibitem{ward}
A.J. Ward, \emph{The topological characterization of an open linear interval}, Proc. London Math. Soc. (2) \textbf{41}
(1936),
191-198.
\end{thebibliography}
\end{document}